\pdfoutput=1
\documentclass[a4paper,reqno,11pt]{amsart}

\usepackage[margin=1in,includehead,includefoot]{geometry}
\usepackage[utf8]{inputenc}
\usepackage[T1]{fontenc}
\usepackage{lmodern,microtype,mathtools,amsthm,amssymb,enumitem,xcolor,tikz,hyperref,bookmark}
\usepackage[nameinlink,capitalise]{cleveref}
\usetikzlibrary{decorations.pathmorphing,decorations.pathreplacing}

\hypersetup{
  pdfauthor={Meike Hatzel, Gwenaël Joret, Piotr Micek, Marcin Pilipczuk, Torsten Ueckerdt, Bartosz Walczak},
  pdftitle={Tight bound on treedepth in terms of pathwidth and longest path},
  colorlinks,
  linkcolor={red!50!black},
  citecolor={blue!50!black},
  urlcolor={blue!80!black}
}

\makeatletter
\let\old@setaddresses\@setaddresses
\def\@setaddresses{\bigskip{\parindent 0pt\let\scshape\relax\let\ttfamily\relax\old@setaddresses}}
\makeatother

\setlist[enumerate]{label=\textup{(\arabic*)},topsep=\medskipamount,noitemsep,leftmargin=*}
\setlist[itemize]{topsep=\medskipamount,noitemsep,labelsep=.6em,leftmargin=1.5em}

\newtheorem{theorem}{Theorem}

\DeclarePairedDelimiter\set{\{}{\}}
\DeclarePairedDelimiter\abs{\lvert}{\rvert}

\DeclareMathOperator{\td}{td}
\DeclareMathOperator{\pw}{pw}
\DeclareMathOperator{\tw}{tw}
\DeclareMathOperator{\depth}{depth}
\DeclareMathOperator{\level}{level}
\DeclareMathOperator{\interior}{int}

\newcommand{\cgB}{\mathcal{B}}
\newcommand{\cgI}{\mathcal{I}}
\newcommand{\bigOh}{\mathcal{O}}
\newcommand{\linkage}{\mathcal{L}}

\let\leq\leqslant
\let\geq\geqslant

\title{Tight bound on treedepth in terms of pathwidth and longest path}

\author[M.~Hatzel\and G.~Joret\and P.~Micek\and M.~Pilipczuk\and T.~Ueckerdt\and B.~Walczak]{Meike Hatzel\and Gwenaël Joret\and Piotr Micek\and Marcin Pilipczuk\and Torsten Ueckerdt\and Bartosz Walczak}

\address[M.~Hatzel]{National Institute of Informatics, Tokyo, Japan}
\email{meikehatzel@nii.ac.jp}
\address[G.~Joret]{Computer Science Department, Université libre de Bruxelles, Brussels, Belgium}
\email{gwenael.joret@ulb.be}
\address[P.~Micek, B.~Walczak]{Department of Theoretical Computer Science, Faculty of Mathematics and Computer Science, Jagiellonian University, Kraków, Poland}
\email{piotr.micek@uj.edu.pl, bartosz.walczak@uj.edu.pl}
\address[M.~Pilipczuk]{Institute of Informatics, University of Warsaw, Warsaw, Poland}
\email{malcin@mimuw.edu.pl}
\address[T.~Ueckerdt]{Computer Science Department, Karlsruhe Institute of Technology, Karlsruhe, Germany}
\email{torsten.ueckerdt@kit.edu}

\thanks{M.~Hatzel was supported by the Federal Ministry of Education and
Research (BMBF) and by a fellowship within the IFI programme of the German Academic Exchange Service (DAAD).
G.~Joret is supported by a CDR grant from the Belgian National Fund for Scientific Research (FNRS), a PDR grant from FNRS, and by the Wallonia Brussels International (WBI) agency.
P.~Micek is supported by the National Science Center of Poland under grant no.\ 2018/31/G/ST1/03718.
B.~Walczak is partially supported by the National Science Center of Poland under grant no.\ 2019/34/E/ST6/00443.
The research of M. Pilipczuk is part of a project that has received funding from the European Research Council (ERC) under the European Union's Horizon 2020 research and innovation programme Grant Agreement 714704.}

\begin{document}

\begin{abstract}
We show that every graph with pathwidth strictly less than $a$ that contains no path on $2^b$ vertices as a subgraph has treedepth at most $10ab$.
The bound is best possible up to a constant factor.
\end{abstract}

\maketitle
\vspace*{-2ex}

\section{Introduction}

Treewidth ($\tw$), pathwidth ($\pw$), and treedepth ($\td$) are among the best-known and most widely studied structural width parameters of graphs.
They are related by the inequalities $\tw(G)+1\leq\pw(G)+1\leq\td(G)$ for every graph $G$.
Moreover, trees have treewidth $1$ and arbitrarily large pathwidth, while paths have pathwidth $1$ and arbitrarily large treedepth.

Treedepth is approximated by the maximum length of a path\footnote{In this paper, we are concerned only about non-induced paths.}: every graph containing an $\ell$-vertex path has treedepth greater than $\log_2\ell$, and every graph with no such path has treedepth less than $\ell$ \cite[Section~6]{NO12}.
Similarly, pathwidth is approximated by the maximum height of a complete binary tree minor: every graph containing a complete binary tree of height $h$ as a minor has pathwidth at least $\lfloor\frac{h}{2}\rfloor$~\cite{Scheffler89}, and every graph with no such minor has pathwidth $\bigOh(2^h)$~\cite{BRST91}.
For both parameters, the exponential gap between the respective lower and upper bounds cannot be avoided, as witnessed by complete graphs.
Treewidth is approximated by the maximum size of a grid minor, but here the gap is polynomial: while every graph containing a $k\times k$ grid as a minor has treewidth at least $k$~\cite{RS86}, every graph with no such minor has treewidth polynomial in $k$~\cite{CC16}.

Kawarabayashi and Rossman~\cite{KR22} showed that treedepth is approximated with polynomial gap by the three above-mentioned obstructions together: every graph with no $k\times k$ grid minor, no height $k$ complete binary tree minor, and no $2^k$-vertex path has treedepth polynomial in $k$.
More specifically, they proved that every graph of treewidth less than $k$ with no height $k$ complete binary tree minor and no $2^k$-vertex path has treedepth $\bigOh(k^5\log^2k)$.
Here are an improvement of this statement and an analogous result relating pathwidth and treewidth:

\begin{theorem}[Czerwiński, Nadara, Pilipczuk~\cite{CNP21}\protect\footnote{In~\cite{CNP21}, the bound is stated in the special case $t=h=b$, but the proof works in general.}]
\label{thm:td_vs_tw}
Every graph of treewidth less than\/ $t$ with no complete binary tree of height\/ $h$ as a minor and no\/ $2^b$-vertex path has treedepth\/ $\bigOh(thb)$.
\end{theorem}

\begin{theorem}[Groenland, Joret, Nadara, Walczak~\cite{GJNW23}]
\label{thm:pw_vs_tw}
Every graph of treewidth less than\/ $t$ with no complete binary tree of height\/ $h$ as a minor has pathwidth\/ $\bigOh(th)$.
\end{theorem}

We complete the picture by proving an analogous result relating treedepth and pathwidth.

\begin{theorem}
\label{thm:main}
Every graph of pathwidth less than\/ $a$ containing no\/ $2^b$-vertex path has treedepth at most\/ $10ab$.
\end{theorem}

Clearly, \cref{thm:pw_vs_tw,thm:main} imply \cref{thm:td_vs_tw}.
On the other hand, \cref{thm:td_vs_tw} implies that every graph of pathwidth less than $a$ containing no $2^b$-vertex path has treedepth $\bigOh(a^2b)$.
This is because every graph with pathwidth less than $a$ has treewidth less than $a$ and contains no complete binary tree of height $2a$ as a minor.
In~\cite{GJNW23}, it was conjectured that the bound on treedepth can be reduced to $\bigOh(ab)$, and \cref{thm:main} provides a proof of this conjecture.

We remark that the bound in \cref{thm:main} is sharp up to a constant factor, which can be seen as follows.
Let $b$ and $c$ be integers with $b>c\geq 1$, and let $a=2^c$.
Consider the graph $G$ obtained from a path on $2^{b-c}$ vertices by replacing each vertex with a clique on $\frac{a}{2}=2^{c-1}$ vertices and replacing each edge by a complete bipartite graph between the two cliques.
Then $\pw(G)=a-1$.
Also, $G$ has $2^{b-1}$ vertices, and thus it has no $2^b$-vertex path.
It can be checked that $G$ has treedepth at least $\frac{a}{2}(b-c)$, which is roughly $ab/2$ when $b\gg c$.
It is shown in~\cite{GJNW23} that the bound in \cref{thm:pw_vs_tw} is also sharp up to a constant factor.
Whether the bound in \cref{thm:td_vs_tw} can be improved remains an open problem.

\section{Preliminaries}

All graphs in this paper are finite and simple, that is, they have no loops or parallel edges.
All logarithms in this paper are to the base $2$.

A \emph{rooted tree} is a tree with one vertex designated as the \emph{root}.
A \emph{rooted forest} is a disjoint union of rooted trees.
We define the \emph{height} of a rooted forest $F$ as the maximum number of vertices on a path from a root to a leaf in $F$.
A vertex $u$ is an \emph{ancestor} of a vertex $v$ in a rooted forest $F$ if $u$ lies on the (unique) path from a root to $v$ in $F$.
A rooted forest $F$ is an \emph{elimination forest} of a graph $G$ if $V(F)=V(G)$ and for every edge $uv$ of $G$, one of the vertices $u$ and $v$ is an ancestor of the other in $F$.
The \emph{treedepth} of a graph $G$, denoted by $\td(G)$, is the minimum height of an elimination forest of $G$.

A \emph{tree decomposition} of a graph $G$ is a pair $(T,\cgB)$ such that $T$ is a tree, the vertices of which are called \emph{nodes}, and $\cgB$ is a collection $\set{B_t}_{t\in V(T)}$ of subsets of $V(G)$, called \emph{bags}, indexed by the nodes of $T$, such that the following conditions are satisfied:
\begin{enumerate}
\item for every edge $uv\in E(G)$, there is a bag containing both $u$ and $v$;
\item for every vertex $v\in V(G)$, the set of nodes $t\in V(T)$ with $v\in B_t$ induces a non-empty subtree of $T$.
\end{enumerate}
The \emph{width} of a tree decomposition is the maximum size of a bag minus $1$.
The \emph{treewidth} of a graph $G$, denoted by $\tw(G)$, is the minimum width of a tree decomposition of $G$.
The notions of \emph{path decomposition} and \emph{pathwidth} are defined analogously with the extra condition that the tree $T$ is a path.
The pathwidth of $G$ is denoted by $\pw(G)$.

A \emph{$k$-linkage} between two subsets $A$ and $B$ of the vertices of a graph $G$ is a subgraph of $G$ that consists of $k$ vertex-disjoint paths each starting in $A$ and ending in $B$.
(If $A$ and $B$ intersect, then a path of a $k$-linkage between $A$ and $B$ may consist of a single vertex in $A\cap B$.)
A path decomposition $(P,\cgB)$ of a graph $G$ is \emph{linked} if for any two nodes $t,t'\in V(P)$, there is a $k$-linkage between $B_t$ and $B_{t'}$ where $k$ is the minimum size of a bag $B_s$ for nodes $s$ on the path from $t$ to $t'$ in $P$.
We use the fact that there is always a path decomposition of minimum width that is linked.

\begin{theorem}[{Erde~\cite[Theorem~5.8]{Erde18}}]
\label{thm:linked}
Every graph\/ $G$ has a path decomposition of width\/ $\pw(G)$ that is linked.
\end{theorem}

\section{Proof}

We proceed with the proof of \cref{thm:main}, that every graph of pathwidth less than $a$ containing no $2^b$-vertex path has treedepth at most $10ab$.

Let $G$ be a graph with $\pw(G)<a$ and with no path on $2^b$ vertices.
If $2^b<2a$, then the statement of the theorem is easily seen to hold by considering a depth-first search forest of $G$, which is an elimination forest of $G$.
Its height is less than $2^b$, which is less than $2a$.
Hence, $\td(G)<2a<10ab$.
Therefore, we may assume that $2^b\geq 2a$.
This inequality will be used at the very end of the proof.

Fix a linked path decomposition $(P,\cgB)$ of $G$ with $\cgB=\set{B_t}_{t\in V(P)}$ and $\abs{B_t}\leq a$ for every node $t\in V(P)$; such a linked path decomposition exists by \cref{thm:linked}.
We think of the nodes as being laid out from left to right along $P$.
For a set of nodes $X\subseteq V(P)$, let
\[B(X) = \bigcup_{t\in X}B_t.\]
We call the node set of any subpath of $P$ an \emph{interval}.
For an interval $I$, we let
\[\level(I) = \min\set{\abs{B_t}\mid t\in I} \qquad \text{and} \qquad \interior(I) = B(I)-B(V(P)-I).\]
Thus, $\interior(I)$ (the ``interior'' of $I$) is the set of vertices of $G$ that lie only in the bags of nodes in $I$.

For every $k\in\set{1,\ldots,a}$ and every inclusion-maximal interval $I^*$ with $\level(I^*)\geq k$, we fix some $k$-linkage between the bags of the leftmost and the rightmost nodes in $I^*$, and we let $\linkage^*_k(I^*)$ be the vertex set of that $k$-linkage.
For every $k\in\set{1,\ldots,a}$ and every interval $I$ with $\level(I)\geq k$, we let $\linkage_k(I)=\linkage^*_k(I^*)\cap B(I)$ where $I^*$ is the unique inclusion-maximal interval with $\level(I^*)\geq k$ containing $I$.
We note the following properties of the sets $\linkage_k(I)$ for further reference:
\begin{gather}
\linkage_k(I)\subseteq B(I),
\label{eq:linkage_in_B} \\
\linkage_k(I')\subseteq\linkage_k(I) \qquad \text{for every interval }I'\subseteq I,
\label{eq:linkage_monotone} \\
\abs{\linkage_k(I)\cap B_t}\geq k \qquad \text{for every node }t\in I,
\label{eq:linkage_wide} \\
\abs{\linkage_k(I)}<k\cdot\smash[t]{2^b}.
\label{eq:linkage_short}
\end{gather}

We describe an iterative process in which we construct a rooted tree $T$ whose vertices are contained in $V(G)$ except for the root, which is a special vertex $r^*\notin V(G)$.
The initial tree $T$ contains only the root $r^*$.
We grow the tree $T$ in rounds, in each round attaching new paths formed by some vertices of $G$ that are not yet in $T$.
We maintain the invariant that $T-r^*$ is an elimination forest of the corresponding induced subgraph of $G$, that is, for any two vertices in $V(T)-\set{r^*}$ that are adjacent in $G$, one is an ancestor of the other in $T$.
The process ends when $T$ contains all vertices of $G$, so that $T-r^*$ is an elimination forest of $G$.

A simple plan for a round would be to find a bag $B_t$ whose removal from $G$ would halve some measure that is proportional to the logarithm of the maximum path length.
Then, after adding $B_t$ to $T$ (as a path), we could continue growing $T$ independently on each of the two sides of $G-B_t$ starting from the vertex of $B_t$ that is currently a leaf of $T$.
This is too simple to work, but it motivates our actual approach.

For a tree $T$ as above and an interval $I$, we use the following notation.
Let $\ell=\level(I)$.
For every $k\in\set{1,\ldots,\ell}$, we define
\[x_k(I,T) = \abs{(\interior(I)\cap\linkage_k(I))-V(T)} \qquad \text{and} \qquad w_k(I,T) = \vphantom{\Big|}\smash[t]{\sum_{i=1}^k\log(x_i(I,T)+1)}.\]
The following ``monotonicity'' property is a direct consequence of \eqref{eq:linkage_monotone}:
\begin{equation}
\text{if $I'\subseteq I$ and $V(T')\supseteq V(T)$, then $x_k(I',T')\leq x_k(I,T)$ and $w_k(I',T')\leq w_k(I,T)$,}
\label{eq:x_w_monotone}
\end{equation}
Furthermore, it follows from \eqref{eq:linkage_short} that $x_i(I,T)+1\leq i\cdot 2^b$ for every $i\in\set{1,\ldots,\ell}$, which yields
\begin{align}
w_\ell(I,T) &\leq \sum_{i=1}^\ell\log(i\cdot 2^b) = b\ell+\log(\ell!),
\label{eq:w_bound} \\
w_\ell(I,T)-w_k(I,T) &\leq \sum_{i=k+1}^\ell\log(i\cdot 2^b) = b(\ell-k)+\log\left(\frac{\ell!}{k!}\right)\quad\text{for every }k\in\set{1,\ldots,\ell}.
\label{eq:w_difference}
\end{align}
For notational convenience, we also define $w_i(\emptyset,T)=0$.

For a vertex $v$ of $G$ that has been added to $T$ at some time in the process, let $\depth(v)$ denote the number of vertices of $G$ on the path from $r^*$ to $v$ in $T$ (thus disregarding $r^*$ in the count).
Since we always augment the tree $T$ by adding new vertices as leaves, $\depth(v)$ is determined when $v$ is added to $T$ and remains unchanged till the end of the process.

During the aforesaid iterative process of constructing the tree $T$, we maintain
\begin{itemize}
\item a set $X$ of nodes of $P$, and the invariant that $B(X)\subseteq V(T)$;
\item the family $\cgI$ of intervals contained in $V(P)-X$ that are inclusion-maximal;
\item a designated vertex $v_I$ in $T$ for every interval $I\in\cgI$.
\end{itemize}
We also maintain the following invariants:
\begin{enumerate}[label=Inv.~\arabic*.,ref=Inv.~\arabic*]
\item\label{inv1} For every interval $I\in\cgI$, the path from the root $r^*$ to $v_I$ in $T$ contains every vertex of $B(I)\cap V(T)$.
\item\label{inv2} For every interval $I\in\cgI$ and for $\ell=\level(I)$, we have
\[\tfrac{1}{5}\depth(v_I)+w_\ell(I,T) \leq (b+1)\ell+\log(\ell!).\]
\end{enumerate}
The former allows us to maintain the invariant that $T-r^*$ is an elimination forest of the corresponding induced subgraph of $G$, and the latter helps us bound the height of $T$.

Initially, the tree $T$ contains only the root $r^*$, the set $X$ is empty, $\cgI=\set{I}$ where $I=V(P)$, and $v_I=r^*$.
For this initial setup, \ref{inv1} holds by the fact that $B(I)\cap V(T)=\emptyset$, whereas \ref{inv2} holds by \eqref{eq:w_bound} and the fact that $\depth(v_I)=0$.

Every round consists in choosing an arbitrary interval $I\in\cgI$ and adding one or two nodes of $I$ to $X$.
As a result, the interval $I$ is replaced in $\cgI$ by at most three of its proper subintervals.

\begin{figure}[p]
\centering
\colorlet{color1}{orange!75!black}
\colorlet{color2}{green!50!black}
\colorlet{color3}{red!75!black}
\colorlet{background1}{orange!15!white}
\colorlet{background2}{green!20!white}
\begin{tikzpicture}[vertex/.style={circle,draw,fill,minimum size=4pt,inner sep=0pt}]
  \fill[background2,rounded corners] (-4.7,-2.6) rectangle (9.7,2.6);
  \fill[background1,rounded corners] (-4.5,-2.4) rectangle (3.5,2.4);
  \draw (-4,0) node[vertex,label=above:$r^*$] {} -- (0,2) -- (0,-2) -- cycle;
  \node at (-1.5,0) {$T$};
  \node[vertex,label=left:$v_I$] (p0) at (0,1) {};
  \begin{scope}[color1]
    \node at (-4,2) {$T'$};
    \node[vertex] (p1) at (1,1) {};
    \node[vertex,label=above:$v'_I$] (p2) at (3,1) {};
    \draw[thick] (p0) -- (p1);
    \draw[thick,decorate,decoration={snake,amplitude=2pt,segment length=3.5mm}] (p1) -- (p2);
    \node[rotate=-70,right] at (p1) {\small\,$\in\linkage_{m+1}(I)-V(T)$};
    \node[rotate=-70,right] at (p2) {\small\,$\in\linkage_\ell(I)-V(T)$};
  \end{scope}
  \begin{scope}[color2]
    \node at (9,-2) {$T''$};
    \node[vertex] (p3) at (4,1) {};
    \node[vertex] (p4) at (6,1) {};
    \node[vertex] (p5) at (7,1) {};
    \node[vertex,label=above:$v''_I$] (p6) at (9,1) {};
    \draw[thick] (p2) -- (p3);
    \draw[thick,decorate,decoration={snake,amplitude=2pt,segment length=3.5mm}] (p3) -- node[below,outer sep=2pt] {$B_{t_1}-V(T')$} (p4);
    \draw[thick] (p4) -- (p5);
    \draw[thick,decorate,decoration={snake,amplitude=2pt,segment length=3.5mm}] (p5) -- node[below,outer sep=2pt] {$B_{t_2}-V(T')$} (p6);
  \end{scope}
\end{tikzpicture}\\[3ex]
\begin{tikzpicture}[xscale=1.2]
  \useasboundingbox (-1.2,-3) rectangle (10.2,2.6);
  \node[right] at (-1.3,2.2) {if $w_\ell(L_t,T')\leq w_\ell(R_t,T')$ for every small node $t\in I$:};
  \draw[thick,color3] (1,0) ellipse (21pt and 50pt) node[yshift=1.2cm] {$\in X$};
  \draw[thick,color2,fill=background2] (4,0) ellipse (21pt and 50pt) +(0,1.2) node {$t_1$};
  \draw[thick,color3] (8,0) ellipse (21pt and 50pt) node[yshift=1.2cm] {$\in X$};
  \foreach\x in {0,2,3,5,6,7,9} \draw (\x,0) ellipse (21pt and 50pt);
  \node[fill=white] at (6,0) {not small};
  \draw[decoration={brace,mirror,amplitude=5pt},decorate] (1.5,-2) -- node[below,outer sep=5pt] {$L=L_{t_1}$} (3.5,-2);
  \draw[decoration={brace,mirror,amplitude=5pt},decorate] (4.5,-2) -- node[below,outer sep=5pt] {$M=R_{t_1}$} (7.5,-2);
\end{tikzpicture}\\[1ex]
\begin{tikzpicture}[xscale=1.2]
  \useasboundingbox (-1.2,-3) rectangle (10.2,2.6);
  \node[right] at (-1.3,2.2) {if $w_\ell(L_t,T')>w_\ell(R_t,T')$ for every small node $t\in I$:};
  \draw[thick,color3] (1,0) ellipse (21pt and 50pt) node[yshift=1.2cm] {$\in X$};
  \draw[thick,color2,fill=background2] (5,0) ellipse (21pt and 50pt) node[yshift=1.2cm] {$t_2$};
  \draw[thick,color3] (8,0) ellipse (21pt and 50pt) node[yshift=1.2cm] {$\in X$};
  \foreach\x in {0,2,3,4,6,7,9} \draw (\x,0) ellipse (21pt and 50pt);
  \node[fill=white] at (3,0) {not small};
  \draw[decoration={brace,mirror,amplitude=5pt},decorate] (1.5,-2) -- node[below,outer sep=5pt] {$M=L_{t_2}$} (4.5,-2);
  \draw[decoration={brace,mirror,amplitude=5pt},decorate] (5.5,-2) -- node[below,outer sep=5pt] {$R=R_{t_2}$} (7.5,-2);
\end{tikzpicture}\\
\begin{tikzpicture}[xscale=1.2]
  \useasboundingbox (-1.2,-3) rectangle (10.2,2.6);
  \node[right] at (-1.3,2.2) {otherwise:};
  \draw[thick,color3] (1,0) ellipse (21pt and 50pt) node[yshift=1.2cm] {$\in X$};
  \draw[thick,color2,fill=background2] (3,0) ellipse (21pt and 50pt) node[yshift=1.2cm] {$t_1$};
  \draw[thick,color2,fill=background2] (6,0) ellipse (21pt and 50pt) node[yshift=1.2cm] {$t_2$};
  \draw[thick,color3] (8,0) ellipse (21pt and 50pt) node[yshift=1.2cm] {$\in X$};
  \foreach\x in {0,2,4,5,7,9} \draw (\x,0) ellipse (21pt and 50pt);
  \node[fill=white] at (4.5,0) {\!not small\!};
  \draw[decoration={brace,mirror,amplitude=5pt},decorate] (1.5,-2) -- node [pos=0.5,below,yshift=-0.2cm] {$L=L_{t_1}$} (2.5,-2);
  \draw[decoration={brace,mirror,amplitude=5pt},decorate] (3.5,-2) -- node [pos=0.5,below,yshift=-0.2cm] {$M=L_{t_2}\cap R_{t_1}$} (5.5,-2);
  \draw[decoration={brace,mirror,amplitude=5pt},decorate] (6.5,-2) -- node [pos=0.5,below,yshift=-0.2cm] {$R=R_{t_2}$} (7.5,-2);
\end{tikzpicture}
\end{figure}

Now, we describe the details of a single round.
Pick an interval $I\in\cgI$.
Let
\[\ell=\level(I) \qquad \text{and} \qquad m=\max(\set{0}\cup\set{i\in\set{1,\ldots,\ell}\mid\linkage_i(I)-V(T)=\emptyset}).\]
It follows that
\begin{equation}
\text{for every node $t\in I$, at least $m$ vertices of $B_t$ are already in $T$.}
\label{eq:m_in_T}
\end{equation}
Indeed, this is true if $m=0$, and if $m>0$, this follows from \eqref{eq:linkage_wide} and the fact that all vertices of $\linkage_m(I)$ are already in $T$.
For every $i\in\set{1,\ldots,\ell}$, since every vertex in $\linkage_i(I)-\interior(I)$ belongs to the bag of a neighbor of $I$ in $P$, which belongs to $X$, we have
\begin{equation}
\linkage_i(I)-V(T)\subseteq\linkage_i(I)-B(X)\subseteq\interior(I).
\label{eq:linkage_internal}
\end{equation}

Let $k=\ell-m$.
(We may have $m=\ell$ and $k=0$.)
For every $i\in\set{m+1,\ldots,\ell}$, by \eqref{eq:linkage_internal}, we have $x_i(I,T)=\abs{\linkage_i(I)-V(T)}\geq 1$.
It follows that
\begin{equation}
w_\ell(I,T) = \sum_{i=1}^\ell\log(x_i(I,T)+1) \geq\! \sum_{i=m+1}^\ell\log(x_i(I,T)+1) \geq \ell-m=k.
\label{eq:w_general_bound}
\end{equation}

Choose one vertex from the set $\linkage_i(I)-V(T)$ for each $i\in\set{m+1,\ldots,\ell}$, and add these at most $k$ vertices into $T$ as a path with one end attached to $v_I$.
That is, the first vertex is added as a child of $v_I$ and every further vertex is added as a child of the previous one.
Let $v'_I$ be the last such vertex (i.e., the other end of the path) if $k\geq 1$, and let $v'_I=v_I$ if $k=0$.
Let $T'$ denote the resulting augmented tree.

Call a node $t\in I$ \emph{small} if $\abs{B_t}\leq\ell+k$.
By the definition of $\level(I)$, at least one node in $I$ is small.
It follows from \eqref{eq:m_in_T} that
\begin{equation}
\text{for each small node $t\in I$, at most $\ell+k-m=2k$ vertices of $B_t$ are not yet in $T'$.}
\label{eq:2k_not_in_T}
\end{equation}
Recall that we think of $I$ as ordered from left to right.
For every node $t\in I$, let $L_t$ and $R_t$ denote the sets of nodes of $I$ to the left and to the right of $t$, respectively, so that $I=L_t\cup\set{t}\cup R_t$.
If $w_\ell(L_t,T')\leq w_\ell(R_t,T')$ for every small node $t\in I$, then let $t_1$ be the rightmost small node in $I$, and let $L=L_{t_1}$ and $M=R_{t_1}$.
Similarly, if $w_\ell(L_t,T')>w_\ell(R_t,T')$ for every small node $t\in I$, then let $t_2$ be the leftmost small node in $I$, and let $M=L_{t_2}$ and $R=R_{t_2}$.
Otherwise, let $t_1$ be the rightmost small node in $I$ such that $w_\ell(L_{t_1},T')\leq w_\ell(R_{t_1},T')$ and $t_2$ be the leftmost small node in $I$ such that $w_\ell(L_{t_2},T')>w_\ell(R_{t_2},T')$.
In this case, by \eqref{eq:x_w_monotone}, $t_1$ and $t_2$ occur in this order from left to right, and there are no small nodes between them.
Now, let $L=L_{t_1}$, $R=R_{t_2}$, and $M=R_{t_1}\cap L_{t_2}$ (i.e., $M$ is the set of nodes strictly between $t_1$ and $t_2$).
See the figure.

Whenever $t_1$ and $t_2$ are defined, we add them to $X$.
We remove $I$ from $\cgI$, and whenever $L$, $M$, and $R$ are defined and non-empty, we add them as new intervals to $\cgI$.

Now, we add the vertices of $B_{t_1}-V(T')$ and $B_{t_2}-V(T')$ (whenever $t_1$ or $t_2$ are defined) to $T'$ as one path with one end attached to $v'_I$.
That is, the first such vertex is a child of $v'_I$, and every further vertex is a child of the previous.
Note that possibly both sets $B_{t_1}-V(T')$ and $B_{t_2}-V(T')$ are empty; in particular, this happens when $k=0$.
Let $v''_I$ be the last vertex added if at least one vertex was added, and let $v''_I=v'_I$ otherwise.
Let $T''$ denote the new tree.
By \eqref{eq:2k_not_in_T}, we have added at most $4k$ additional vertices, so
\begin{equation}
\depth(v''_I) \leq \depth(v'_I)+4k = \depth(v_I)+5k.
\label{eq:inequality_on_depth}
\end{equation}
Whenever $L$, $M$, or $R$ is defined and non-empty, we set the corresponding vertex $v_L$, $v_M$, or $v_R$ to be $v''_I$.
By \ref{inv1} for $I$, it follows that the path from $r^*$ to $v''_I$ contains every vertex of $B(I)\cap V(T'')$, which yields \ref{inv1} for $L$, $M$, and $R$ (when they are defined and non-empty).

Before verifying \ref{inv2}, let us capture the key properties of $L$, $M$, and $R$.
If $L$ is defined and non-empty (so that $L=L_{t_1}$), then let $\overline{L}=R_{t_1}$.
If $R$ is defined and non-empty (so that $R=R_{t_2}$), then let $\overline{R}=L_{t_2}$.
Whenever the respective sets are defined, we have
\begin{gather}
w_\ell(L,T') \leq w_\ell(\overline{L},T') \qquad\text{and}\qquad w_\ell(R,T') \leq w_\ell(\overline{R},T'), \label{eq:weight_of_L_R_aux}\\
w_\ell(M,\smash[t]{T''}) \leq w_\ell(I,T), \label{eq:weight_of_M}\\
\level(M) \geq \ell+k+1, \label{eq:level_of_M}
\end{gather}
where \eqref{eq:weight_of_M} follows from \eqref{eq:x_w_monotone}, and \eqref{eq:level_of_M} follows as there are no small nodes in $M$.

While a bound analogous to \eqref{eq:weight_of_M} holds also for $w_\ell(L,T'')$ and $w_\ell(R,T'')$, we need a stronger one.
First we focus on the interval $L$, and the argument is symmetric for the interval $R$.
For $i\in\set{1,\ldots,\ell}$, we compare $x_i(I,T)$ with $x_i(L,T')$ and $x_i(\overline{L},T')$.
Note that $\interior(L)$ and $\interior(\overline{L})$ are vertex-disjoint and are both contained in $\interior(I)$.
For each $i\in\set{1,\ldots,\ell}$, we have $\linkage_i(L)\subseteq\linkage_i(I)$ and $\linkage_i(\overline{L})\subseteq\linkage_i(I)$, by \eqref{eq:linkage_monotone}.
For each $i\in\set{m+1,\ldots,\ell}$, we have put one vertex of $\linkage_i(I)-V(T)$ into $T'$; this vertex belongs to $\interior(I)$ by \eqref{eq:linkage_internal}.
This, the property \eqref{eq:linkage_in_B}, and the definition of $x_i$ imply that for each $i\in\set{1,\ldots,\ell}$, we have
\[x_i(I,T) \geq \begin{cases}
x_i(L,T')+x_i(\overline{L},T') & \text{if }i\leq m,\\
x_i(L,T')+x_i(\overline{L},T')+1 & \text{if }i>m.
\end{cases}\]
Since $x_i(L,T')$ and $x_i(\overline{L},T')$ are non-negative, the above implies
\begin{equation}
x_i(I,T)+1 \geq \begin{cases}
x_i(L,T')+x_i(\overline{L},T')+1 \geq \tfrac{1}{2}(x_i(L,T')+x_i(\overline{L},T')+2) & \text{if }i\leq m,\\
x_i(L,T')+x_i(\overline{L},T')+2 & \text{if }i>m.
\end{cases}
\label{eq:x_lower_bound}
\end{equation}
Recalling that $k=\ell-m$, we calculate
\begin{alignat*}{2}
w_\ell(I,T) &= \sum_{i=1}^\ell\log(x_i(I,T)+1)\\
&\geq \sum_{i=1}^\ell\log(x_i(L,T')+x_i(\overline{L},T')+2)-m && \quad\text{by \eqref{eq:x_lower_bound}}\\
&\geq \sum_{i=1}^\ell\frac{\log(x_i(L,T')+1)+\log(x_i(\overline{L},T')+1)}{2}+\ell-m && \quad\text{($*$)}\\
&= \tfrac{1}{2}(w_\ell(L,T')+w_\ell(\overline{L},T'))+k\\
&\geq w_\ell(L,T')+k && \quad\text{by \eqref{eq:weight_of_L_R_aux}}\\
&\geq w_\ell(L,T'')+k && \quad\text{by \eqref{eq:x_w_monotone}},
\end{alignat*}
where in ($*$), we use the inequality $\log(x+y)=\log\frac{x+y}{2}+1\geq\frac{1}{2}(\log x+\log y)+1$ that follows from the concavity of $\log$.
From this and the analogous argument for $R$, we conclude that
\begin{equation}
w_\ell(L,T'')+k\leq w_\ell(I,T) \qquad\text{and}\qquad w_\ell(R,T'')+k\leq w_\ell(I,T).
\label{eq:weight_of_L_R}
\end{equation}

Now, we are set to verify \ref{inv2} for the intervals $L$, $R$, and $M$ (when they are defined and non-empty).
We have
\begin{alignat*}{2}
&\tfrac{1}{5}\depth(v_L)+w_{\level(L)}(L,T'')\\[-.5ex]
&\quad\leq \tfrac{1}{5}\depth(v_I)+k+w_\ell(L,T'')+b(\level(L)-\ell)+\log\left(\frac{\level(L)!}{\ell!}\right) && \quad\text{by \eqref{eq:inequality_on_depth} and \eqref{eq:w_difference}}\\
&\quad\leq \tfrac{1}{5}\depth(v_I)+w_\ell(I,T)+b(\level(L)-\ell)+\log\left(\frac{\level(L)!}{\ell!}\right) && \quad\text{by \eqref{eq:weight_of_L_R}}\\
&\quad\leq (b+1)\ell+\log(\ell!)+b(\level(L)-\ell)+\log\left(\frac{\level(L)!}{\ell!}\right) && \quad\text{by \ref{inv2} for }I\\
&\quad\leq (b+1)\level(L)+\log(\level(L)!).
\end{alignat*}
The exact same bounds hold with $L$ replaced by $R$.
Finally, for $M$, we have
\begin{alignat*}{2}
&\tfrac{1}{5}\depth(v_M)+w_{\level(M)}(M,T'')\\[-.5ex]
&\quad\leq \tfrac{1}{5}\depth(v_I)+k+w_\ell(M,T'')+b(\level(M)-\ell)+\log\left(\frac{\level(M)!}{\ell!}\right) && \quad\text{by \eqref{eq:inequality_on_depth} and \eqref{eq:w_difference}}\\
&\quad\leq \tfrac{1}{5}\depth(v_I)+k+w_\ell(I,T)+b(\level(M)-\ell)+\log\left(\frac{\level(M)!}{\ell!}\right) && \quad\text{by \eqref{eq:weight_of_M}}\\
&\quad\leq (b+1)\ell+\log(\ell!)+k+b(\level(M)-\ell)+\log\left(\frac{\level(M)!}{\ell!}\right) && \quad\text{by \ref{inv2} for }I\\
&\quad\leq (b+1)\ell+\level(M)-\ell+b(\level(M)-\ell)+\log(\level(M)!) && \quad\text{by \eqref{eq:level_of_M}}\\
&\quad= (b+1)\level(M)+\log(\level(M)!).
\end{alignat*}

This completes the round of our process for the interval $I$, with $T''$ becoming the new tree $T$.
We have shown that both invariants, \ref{inv1} and \ref{inv2}, are preserved.

The process ends when all vertices of $G$ have been added to $T$.
It remains to show that $T-r^*$ is an elimination forest of $G$ with height at most $10ab$.
To see that it is an elimination forest, observe that whenever a vertex $v\in\interior(I)$ is added to $T$ when considering an interval $I$, \ref{inv1} guarantees that all neighbors of $v$ in $G$ that are already in $T$ lie on the path from $r^*$ to $v$ in $T$, as the neighbors of $v$ in $G$ belong to $B(I)$ by the definition of path decomposition.

The height of the forest $T-r^*$ is equal to $\max\set{\depth(v)\mid v\in V(G)}$.
Let $v$ be a vertex of $G$.
Consider the moment in the process when $v$ has been added to $T$.
Say, it happened when processing an interval $I$ with $\level(I)=\ell$.
Let $m$ and $k$ be the values fixed when processing $I$.
Clearly, we have $\depth(v)\leq\depth(v''_I)$ and $\ell\leq a$.
Therefore,
\begin{alignat*}{2}
\depth(v) &\leq \depth(v''_I)\\
&\leq \depth(v_I)+5k && \quad\text{by \eqref{eq:inequality_on_depth}}\\
&\leq \depth(v_I)+5w_\ell(I,T) && \quad\text{by \eqref{eq:w_general_bound}}\\
&\leq 5(b+1)\ell+5\log(\ell!) && \quad\text{by \ref{inv2}}\\
&\leq 5(b+1)a+5\log(a!)\\
&\leq 5ab+5a\log a+5a.
\end{alignat*}

Recall that $2^b\geq 2a$ and thus $b\geq\log(2a)=\log a+1$.
It follows that
\[\depth(v)\leq 5ab+5a\log a+5a\leq 10ab.\]
We conclude that $T-r^*$ is an elimination forest of $G$ with height at most $10ab$, as desired.

\pagebreak

\section*{Acknowledgements}

This research was carried out at the ``Sparsity in Algorithms, Combinatorics, and Logic'' workshop held in Dagstuhl in September 2021.
We thank the organizers and the workshop participants for creating a productive working atmosphere.
In particular, we thank Marthe Bonamy, Marcin Briański, Kord Eickmeyer, Wojciech Nadara, Ben Rossman, Blair~D. Sullivan, and Alexandre Vigny for initial discussions on the topic of the paper.
We also thank an anonymous reviewer for noticing an error in a previous version of the paper.

\bibliographystyle{plain}
\bibliography{bibliography}

\end{document}